\documentclass[12pt]{amsart}
\usepackage{amsfonts, amsmath, amsthm}

\newcommand\AAA{\mathbb{A}}

\newcommand\PP{\mathbb{P}}

\newcommand\calL{\mathcal{L}}

\newtheorem{theorem}{Theorem}
\newtheorem{lemma}[theorem]{Lemma}

\begin{document}

\title[More \'etale covers of affine spaces]{More \'etale 
covers of affine spaces in positive characteristic}
\author{Kiran S. Kedlaya}
\address{Department of Mathematics, Room 2-165 \\
Massachusetts Institute of Technology \\
77 Massachusetts Avenue \\
Cambridge, MA 02139}
\email{kedlaya@mit.edu}
\urladdr{www-math.mit.edu/\~{}kedlaya}
\thanks{Supported by a National Science Foundation postdoctoral
fellowship. Thanks to Bas Edixhoven and Bjorn Poonen for helpful discussions.}

\date{December 14, 2003}

\begin{abstract}
We prove that every geometrically reduced
projective variety of pure dimension $n$
over a field of positive characteristic
admits a morphism to projective $n$-space, \'etale away from the hyperplane
$H$ at infinity, which maps a chosen divisor into $H$
and some chosen smooth points not on the divisor to points
not in $H$. This improves an earlier result of the author, which
was restricted to infinite perfect fields. We also prove a related result
that controls the behavior of divisors through the chosen point.
\end{abstract}

\maketitle

\section{Results}

We prove the following theorem, which improves upon the main result
of \cite{meetale} by elimining the hypothesis that $k$ be infinite
and perfect. (The case $n=1$ seems to have become a folk theorem,
but we believe it is originally due to Abhyankar.)
\begin{theorem} \label{thm:main}
Let $X$ be a geometrically reduced,
projective variety of pure dimension $n$ over a field $k$ of
characteristic $p>0$. Let $\calL$ be an ample line bundle of
$X$, let $D$ be a closed subscheme of $X$ of dimension less than $n$,
and let $S$ be a
zero-dimensional subscheme of the smooth locus of $X$ not meeting $D$.
Then there exists a finite morphism
$f: X \to \PP^n_k$ of $k$-schemes satisfying the following conditions:
\begin{enumerate}
\item[(a)]
$f$ is induced by an $(n+1)$-tuple of sections of some tensor power of
$\calL$;
\item[(b)]
$f$ is \'etale away from the hyperplane $H \subseteq \PP^n$ at infinity;
\item[(c)]
$f(D) \subseteq H$;
\item[(d)]
$f(S)$ does not meet $H$.
\end{enumerate}
\end{theorem}
Note that the hypothesis that $X$ is geometrically reduced is needed for 
each irreducible component of $X$ to be generically smooth, as any
scheme admitting a generically \'etale map to a smooth scheme must be.

Theorem~\ref{thm:main}
implies in particular that the smooth locus of
$X$ is covered by open affines which admit
finite \'etale morphisms to $\AAA^n$. We can also prove a refinement
of this corollary.
\begin{theorem} \label{thm:main2}
Let $X$ be a separated scheme of finite type
over $k$ (still of characteristic $p>0$)
 of pure dimension $n$, let $x$ be a smooth point
of $X$, and let $D_1, \dots, D_m$ be irreducible divisors of $X$
intersecting transversely at $x$.
Then there exists a finite \'etale morphism $f: U \to \AAA^n$
for $U$ some open dense subset of $X$, defined over $k$ and containing $x$,
such that $D_1, \dots, D_m$ map to coordinate hyperplanes.
\end{theorem}
Note that one cannot expect $f$ to extend to a finite morphism from $X$
to $\PP^n$, especially at points where the $D_i$ have nontransverse
intersections.

\section{Review of Noether normalization}

Before proving the theorems, we review Noether normalization and prove
a strong form of it for our purposes. Most of these results are evident
when $k$ is infinite; the care is required in handling $k$ finite, since
``generic'' constructions are not available.

We first point out a simple fact that will come up repeatedly in what
follows: if $X$ is a projective variety
and $\calL$ is an ample line bundle on $X$, then the complement of the zero
locus of any section of $\calL$ is an affine scheme.

\begin{lemma} \label{lem:noeth}
Let $X$ be a projective scheme over
a field $k$ of dimension $n$, and let
$\calL$ be an ample line bundle on $X$.
Suppose $\calL$ admits linearly independent
sections $\alpha_0, \dots, \alpha_n$ whose
zero loci have no common intersection. Then the $\alpha_i$ induce
a finite map $X \to \PP^n$.
\end{lemma}
\begin{proof}
The $\alpha_i$ clearly induce a map $X \to \PP^n$; all that needs to
be checked is that the map is finite.
Let $t_0, \dots, t_n$ be homogeneous coordinates on $\PP^n$ which
pull back to $\alpha_0, \dots, \alpha_n$ on $X$. Given a geometric point
$x$ of $\PP^n$, choose $i$ such that $t_i$ is nonzero at $x$. Then on
one hand, the fibre of $x$ is projective over $x$, since it is the common
zero locus of certain sections of $\calL$ and so is closed in $X$. 
On the other hand, the fibre of $x$ is affine over $x$, since it
is contained in the affine subset
of $X$ on which $t_i$ does not vanish.
 Hence the fibre of $x$ is finite over $x$, i.e.,
$f$ is quasi-finite. Moreover, $f$ is also proper (because any map between
proper $k$-schemes is itself proper), so $f$ is finite.
\end{proof}

In the next few lemmas, we adopt the following convention.
We say that a statement about a positive integer $n$ is
true ``for $n$ sufficiently divisible'' if it holds whenever $n$ is divisible
by
some (unspecified) positive integer. In particular, any statement that holds
``for $n$ sufficiently large'' also holds ``for $n$ sufficiently divisible''.

\begin{lemma} \label{lem:sep}
Let $X$ be a projective scheme over
a field $k$, and
let $\calL$ be an ample line bundle on
$X$.
Let $D$ be a closed subscheme of $X$, and
let $Z$ be a zero-dimensional closed subscheme of $X$ not meeting $D$.
Then for $l$ sufficiently divisible, there is a section of 
$\calL^{\otimes l}$ vanishing along $D$ but not vanishing at any point
of $Z$.
\end{lemma}
\begin{proof}
For any sufficiently large
$a$, $\calL^{\otimes a}$ admits a section $s_0$ not vanishing
at any point of $Z$. (Namely, this occurs when $\calL^{\otimes a}$
is very ample and has at least $\deg(Z) + 1$ linearly independent sections.)
Since the complement of the zero locus of $s_0$ is affine, it admits a
regular function vanishing along $D$ but not on $S$. That function has
the form $s_1 s_0^{-b}$ for some positive integer $b$; 
then for $l$ divisible by
$ab$, the section $s_1^{\otimes l/(ab)}$ of $\calL^{\otimes l}$
has the desired property.
\end{proof}

\begin{lemma} \label{lem:ample}
Let $X$ be a projective scheme over
a field $k$ of dimension $n$, 
and let $\calL$ be an ample line bundle on $X$.
Let $S$ be a zero-dimensional subscheme
of $X$, and
let $D$ be a closed subscheme of $X$ not meeting $S$,
Suppose that $0 \leq m \leq n$ and that
$D_1, \dots, D_m$ are divisors on $X$, such that for any nonempty subset
$T$ of $\{1, \dots, m\}$,
the intersection 
$D \cap \bigcap_{t \in T} D_t$ has codimension in $D$
at least $\# T$.
Then for $l$ sufficiently divisible,
there exist sections $s_1, \dots, s_n$
of $\calL^{\otimes l}$ with no common zero on $D$, such that
each $s_i$ vanishes on $S$, and $s_i$ vanishes along $D_i$ for $i=1, \dots, 
m$.
\end{lemma}
Here we take the codimension of the empty scheme
in any other scheme to be $+\infty$.
\begin{proof}
We construct the desired sections inductively as follows.
Suppose that for $0 \leq j < n$,
there exist positive integers $l_i$ and sections $s'_i$
of $\calL^{\otimes l_i}$ for $i=1, \dots, j$ satisfying the following 
conditions.
\begin{enumerate}
\item[(a)] Each $s'_i$ vanishes on $S$.
\item[(b)] If $i \leq \min\{j,m\}$, then $s'_i$ vanishes along $D_i$.
\item[(c)] For any subset $T$ of $\{j+1, \dots, m\}$
(which is the empty set if $j > m$),
the intersection $Y_{j,T}$ of $D$, the zero loci of $s'_1, \dots, s'_j$,
and the $D_t$ for $t \in T$, has codimension in $D$ at least $j + \# T$.
\end{enumerate}
For instance, this is true by hypothesis for $j=0$.

With the $s'_i$ as above,
let $Z_j$ be a zero-dimensional subscheme of $D \setminus D_{j+1}$
meeting each irreducible component of $Y_{j,T}$ having codimension $j+ \# T$ in
$D$, for each subset $T$ of $\{j+2, \dots, m\}$.
(By (c), none of these components is actually
contained in $D \cap D_{j+1}$, so a suitable $Z_j$ can be found.)
For $l_{j+1}$ sufficiently divisible, we can then find by 
Lemma~\ref{lem:sep} a section
$s'_{j+1}$ of $\calL^{\otimes l_{j+1}}$ vanishing along $S \cup D_{j+1}$
but not vanishing at any point of $Z_j$. Then conditions (a), (b), (c)
are satisfied with $j$ replaced by $j+1$.

By induction, we can satisfy (a), (b), (c) with $j=n$.
Let $m$ be the least common multiple of the $l_i$; then for any $l$
divisible by $m$,
the sections $s_i = (s'_i)^{l/l_i}$ of $\calL^{\otimes l}$ have the desired
properties.
\end{proof}

\begin{lemma} \label{lem:sections}
Let $X$ be a geometrically reduced projective scheme over
a field $k$ of pure dimension $n$, 
and let $\calL$ be an ample line bundle on $X$.
Let $S$ be a zero-dimensional subscheme of the smooth locus of $X$,
and let $\alpha$ be a section of $\calL$ whose zero locus $D$
does not meet $S$.
Then for $l$ sufficiently divisible,
there exist
sections $\delta_1, \dots, \delta_n$ of $\calL^{\otimes l}$ 
satisfying the following conditions.
\begin{enumerate}
\item[(a)] The $\delta_i$ have no common zero on $D$,
so that $\alpha^l, \delta_1, \dots, \delta_n$ define a finite morphism
$f: X \to \PP^n$ by Lemma~\ref{lem:noeth}.
\item[(b)] The map $f$ is unramified at each point of $S$.
\end{enumerate}
If moreover $D_1, \dots, D_m$ are divisors as in Lemma~\ref{lem:ample},
which additionally meet transversely at each point of $S$,
then 
we can also ensure that $\delta_i$ vanishes along $D_i$
for $i=1, \dots, m$.
\end{lemma}
\begin{proof}
Since $\calL$ is ample, $U = X \setminus D$ is affine.
Since $X$ is smooth and the $D_i$
meet transversely at each point of 
$S$, we can find regular functions $f_1, \dots, f_n$ on $U$, 
with $f_i$ vanishing
along $D_i$ for $i=1, \dots, m$, which induce a map
$U \to \AAA^n$ unramified at each point of $S$. 
(Besides the vanishing conditions along the $D_i$, we need only ensure that
$f_1, \dots, f_n$ span the cotangent space of $X$ at each point of $S$.)
Write each $f_i$ as $\beta_i/\alpha^{j_i}$ for some positive integer
$j_i$ and some section $\beta_i$ of $\calL^{\otimes j_i}$.

By Lemma~\ref{lem:ample}, we can
choose sections $\gamma_1, \dots, \gamma_n$ of $\calL^{\otimes l}$, for $l$
sufficiently divisible, such that each $\gamma_i$ vanishes on $S$,
$\gamma_i$ vanishes along $D_i$ for
$i=1, \dots, m$,
and $\gamma_1, \dots, \gamma_n$ have no common zero
on $D$.
Now put
\[
\delta_i = \beta_i \alpha^{l-2j_i} + \gamma_i^2 \qquad (i=1, \dots, n).
\]
Then the $\delta_i$ have the desired properties.
\end{proof}

\section{Proofs of the theorems}

We now proceed to the proofs of the theorems.
As noted above, Theorem~\ref{thm:main} was established for $k$ infinite
and perfect in \cite{meetale}. However, the proof below is not only more
widely applicable but also somewhat simpler.
This is because we take advantage of a more flexible construction.
The basic idea is that adding a $p$-th power to a function does not
change its ramification on the source, but does move things around 
on the target.

\begin{proof}[Proof of Theorem~\ref{thm:main}]
We may enlarge $S$ so that it meets each irreducible component of $X$.
By replacing $\calL$ by a suitable tensor power, we may also assume 
(thanks to Lemma~\ref{lem:sep}) that
$\calL$ admits a section $s$ vanishing on $D$ but not at any point
of $S$.
For some integer $m$, by Lemma~\ref{lem:sections}
we can find sections $s_1, \dots, s_n$ of $\calL^{\otimes m}$ such that:
\begin{enumerate}
\item[(a)]
The sections $s_1, \dots, s_n$ have no common zero
on $D$, so that $s^m, s_1, \dots, s_n$
induce a finite morphism $g: X \to \PP^n$ by Lemma~\ref{lem:noeth}.
\item[(b)]
The map $g$ is unramified at each point of $S$.
\end{enumerate}
The locus on $X \setminus D$ where $g$ is unramified
is open, and its intersection with each irreducible component of $X$ is
nonempty (since $g$ is unramified on $S$); let $E$ be its complement in $X$.
For some positive integer $r$, by Lemma~\ref{lem:sep} we can find
a section $t$ of $\calL^{\otimes rm}$ which vanishes on $E$ but
not at any point of $S$. 
By Lemma~\ref{lem:ample}, we can find a positive
integer $l \geq 2$ and
sections $t_1, \dots, t_n$ of $\calL^{\otimes lrm}$ which have
no common zero on the vanishing locus $Z$ of
$t$. Now put
\begin{align*}
u_0 &= t^{pl} \\
u_i &= s_i s^{m(pr-1)} t^{p(l-1)} + t_i^p \qquad (i=1, \dots, n).
\end{align*}
Then $u_0, u_1, \dots, u_n$ have no common zero, so they define a 
finite morphism $f: X \to \PP^n$ by Lemma~\ref{lem:noeth}.

For a point $y \in X \setminus Z$, the map $g$ is unramified at
$y$; 
in other words, the differentials of the functions $s_1/s^m, \dots, 
s_n/s^m$ at $y$
are linearly independent. But after multiplication by
$(s^{rm}/t)^p$, these become precisely the
differentials
of the functions $u_1/u_0, \dots, u_n/u_0$. Hence $f$ is unramified
at $y$. That is, $f$ is \'etale over the complement of the hyperplane
$u_0 = 0$, and $D$ maps into this hyperplane but no point of $S$ does.
Thus we have constructed the desired map.
\end{proof}

The proof of Theorem~\ref{thm:main2} is similar, but we must do some
blowing up first.
\begin{proof}[Proof of Theorem~\ref{thm:main2}]
Since we can replace $X$ by a
projective compactification of an affine open neighborhood of $x$,
it suffices to consider $X$ irreducible and projective.
Moreover, by blowing up away from $x$, we can ensure that the intersection
of any subset of the $D_i$ is irreducible.

Choose an ample line bundle $\calL$ on $X$.
We may now repeat the argument of Theorem~\ref{thm:main}
with $S = \{x\}$, while imposing the additional
restriction that $s_i$ and $t_i$ both vanish along
$D_i$. The condition of Lemma~\ref{lem:ample} holds because
for any subset $T \subset \{1, \dots, m\}$,
the intersection $\cap_{t \in T} D_t$ has been arranged to be
irreducible, and does not lie in $D$ or $Z$ because it
contains $x$.
\end{proof}

\end{document}